\documentclass[reqno,11pt]{amsart}
\usepackage{amsmath,amsrefs,amssymb}

\numberwithin{equation}{section}

\DeclareMathOperator{\Real}{Re}
\DeclareMathOperator{\Imag}{Im}

\DeclareMathOperator{\smallo}{o}
\newcommand{\R}{\mathbb{R}}
\newcommand{\C}{\mathbb{C}}
\renewcommand{\H}{\mathbb{H}}

\newcommand{\N}{\mathbb{N}}

\newcommand{\<}{\left<}
\renewcommand{\>}{\right>}
\renewcommand{\[}{\left[}
\renewcommand{\]}{\right]}
\renewcommand{\(}{\left(}
\renewcommand{\)}{\right)}
\usepackage{xcolor}

\newtheorem{theorem}{Theorem}[section]
\newtheorem{lemma}[theorem]{Lemma}
\newtheorem{remark}[theorem]{Remark}

\begin{document}

\title[Liouville-type results for the CR Yamabe equation in $\H^n$]{Liouville-type results for the CR Yamabe equation in the Heisenberg group}

\author{Joshua Flynn}

\address{Joshua Flynn, Department of Mathematics and Statistics, McGill University, 805 Sherbrooke Street West, Montreal, Quebec H3A 0B9, Canada}
\email{joshua.flynn@mcgill.ca}

\author{J\'er\^ome V\'etois}

\address{J\'er\^ome V\'etois, Department of Mathematics and Statistics, McGill University, 805 Sherbrooke Street West, Montreal, Quebec H3A 0B9, Canada}
\email{jerome.vetois@mcgill.ca}

\thanks{The second author was supported by the NSERC Discovery Grant RGPIN-2022-04213.}

\date{November 10, 2023}

\begin{abstract}
We obtain Liouville-type results for solutions to the CR Yamabe equation in $\H^n$, which extend a result obtained by Jerison and Lee for solutions in $L^{2+2/n}\(\H^n\)$. We obtain our results under either pointwise conditions or integral conditions at infinity. In particular, our results hold for all bounded solutions when $n=2$ and solutions satisfying a pointwise decay assumption when $n\ge3$. The proofs rely on integral estimates combined with a suitable divergence formula.
\end{abstract}

\maketitle

\section{Introduction and main result}\label{Sec1}

We consider the CR Yamabe equation 
\begin{equation}\label{Th1Eq1}
4\Delta_bu=n^2u^{\frac{n+2}{n}}\quad\text{in }\H^n,
\end{equation}
where $\H^n = \mathbb{C}^{n} \times \R$ is the Heisenberg group and $\Delta_bu:=-\Real u_{\alpha\overline\alpha}$ is the Heisenberg sub-Laplacian (definitions are recalled in Section~\ref{Sec2}). Jerison and Lee~\cites{JL2} obtained that the only solutions $u\in L^{2+2/n}\(\H^n\)$ to \eqref{Th1Eq1} are the functions 
\begin{equation}\label{Th1Eq2}
u\(z,t\)=\frac{\(4\Imag\lambda-\left|\mu\right|^2\)^{n/2}}{\left|t+\sqrt{-1}\left|z\right|^2+\<\mu,z\>+\lambda\right|^n}\quad\forall\(z,t\)\in\H^n,
\end{equation}
where $\mu\in\C^n$ and $\lambda\in\C$ are such that $\left|\mu\right|^2<4\Imag\lambda$. Moreover, these functions correspond to the extremal functions for the Sobolev inequality in the Heisenberg group. Garofalo and Vassilev~\cite{GV} extended Jerison and Lee's result to groups of Heisenberg-type and partially symmetric solutions, namely solutions which are invariant with respect to the action of the orthogonal group in the horizontal layer of the Lie algebra. Due to the lack of invariance of the sub-Laplacian with respect to the standard reflections about hyperplanes, works in the literature have only been able to use moving plane arguments in restricted cases in this setting (see for example Birindelli and Prajapat~\cite{BP} and Garofalo and Vassilev~\cite{GV}). This is in contrast with the situation for the Yamabe equation in the Euclidean space (see Caffarelli, Gidas and Spruck~\cite{CGS} and Chen and Li~\cite{CLi}).

\smallskip
In this article, we extend Jerison and Lee's result to solutions satisfying either pointwise conditions or integral conditions at infinity. Our first result is the following:

\begin{theorem}\label{Th1}
Let $n\ge2$ and $u$ be a positive solution to \eqref{Th1Eq1} such that
\begin{equation}\label{Th1Eq3}
u\(z,t\)\le C\(\left|z\right|^2+\left|t\right|\)^{-\frac{n-2}{2}}\quad\forall\(z,t\)\in\H^n\backslash\left\{\(0,0\)\right\}
\end{equation}
for some constant $C>0$. Then $u$ is of the form \eqref{Th1Eq2}.
\end{theorem}

\begin{remark}
Catino, Li, Monticelli and Roncoroni~\cite{CLMR} recently posted an article on arXiv, where they obtained a result in the same direction as Theorem~\ref{Th1}. In the case where $n=1$, their result provides the remarkable full classification of solutions to \eqref{Th1Eq1}. In the case where $n\ge2$, they obtain the result under the assumption that 
$$u\(z,t\)\le C\(\left|z\right|^2+\left|t\right|\)^{-\frac{n}{2}}\quad\forall\(z,t\)\in\H^n\backslash\left\{\(0,0\)\right\}$$
for some constant $C>0$. In a first version of our paper, independently from~\cite{CLMR}, we obtained Theorem~\ref{Th1} under the assumption that 
$$u\(z,t\)\le C\(\left|z\right|^2+\left|t\right|\)^{-p}\quad\forall\(z,t\)\in\H^n\backslash\left\{\(0,0\)\right\}$$
for some constant $C>0$ and $p>\frac{n\(n-2\)}{2\(n-1\)}$. This can be improved to $p\ge\frac{n-2}{2}$. We point out that to obtain the limit case where $p=\frac{n-2}{2}$, we now use an argument from~\cite{CLMR} in Lemma~\ref{Lem2}. 
\end{remark}

\begin{remark}
The classification of bounded solutions for $n=2$ in Theorem~\ref{Th1} is potentially useful to obtain a priori estimates near isolated blow-up points for solutions to Yamabe-type equations in CR manifolds of dimension five, in a similar way as in the Riemannian setting (see for instance Schoen~\cites{S1,S2}, Li~\cites{Li1,Li2}, Chen and Lin~\cites{CLin1,CLin2}, Li and Zhu~\cite{LiZhu}, Druet~\cite{D}, Marques~\cite{Ma}, Li and Zhang~\cites{LiZha1,LiZha2,LiZha3} and Khuri, Marques and Schoen~\cite{KMS}). 
\end{remark}

\smallskip
Our second result (of which Theorem~\ref{Th1} is in fact a corollary) is the following:

\begin{theorem}\label{Th2}
Let $n\ge2$ and $u$ be a positive solution to \eqref{Th1Eq1} such that
\begin{equation}\label{Th2Eq1}
\int_{B_R\(0\)}u^{q}\le CR^2\quad\forall R>1
\end{equation}
for some constants $C>0$ and $q\in\(\frac{2n+1}{n},\frac{2n+2}{n}\]$.\,Then $u$ is of the form \eqref{Th1Eq2}.
\end{theorem}

It is interesting to compare \eqref{Th2Eq1} with the estimate
\begin{equation}\label{Th2Eq2}
\int_{B_R\(0\)}u^{q}\le CR^{2n+2-nq}\quad\forall R>1,
\end{equation}
which holds for all $q\in\[0,\frac{n+2}{n}\]$ and all positive solutions to \eqref{Th1Eq1}, without any further assumptions (see Lemma~\ref{Lem5}). Theorem~\ref{Th1} follows from Theorem~\ref{Th2} by showing that the estimate \eqref{Th2Eq1} with $q=\frac{2n+2}{n}$ can be obtained by putting together the pointwise estimate \eqref{Th1Eq3} with the integral estimate \eqref{Th2Eq2}. Theorem~\ref{Th2} is more general than Theorem~\ref{Th1}: it applies to smaller values of $q$ and solutions which do not satisfy a pointwise decay condition.

\smallskip
We prove Theorems~\ref{Th1} and~\ref{Th2} in Section~\ref{Sec2}. The proofs use an approach based on integral estimates. An approach of this type has recently been used by Ma and Ou~\cite{MO} for equations of type (1.1) with subcritical exponents in $\H^n$ and Catino, Monticelli and Roncoroni~\cite{CMR} and Ou~\cite{Ou} (see also V\'etois~\cite{Vet}) for the critical $p$-Laplace equation in $\R^n$. The starting point of this method is to obtain a suitable divergence formula (see Lemma~\ref{Lem1}). In our case, this formula is derived from and extends a celebrated formula introduced and used by Jerison and Lee~\cite{JL2} to characterize the extremal functions of Sobolev inequalities in the Heisenberg group and the CR Yamabe problem. Our extended version of this formula is comparable with the remarkable Obata-type formula found by Ou~\cite{Ou}*{Proposition~2.3} in the context of the critical $p$--Laplace equation in $\R^n$ (see also~\cite{Ob} for the original formula discovered by Obata in his work on the conformal transformations of the sphere). Once our divergence formula is established, we multiply it by some cutoff functions, integrate in $\H^n$ and estimate some remainder terms in a series of technical lemmas (see Lemmas~\ref{Lem2},~\ref{Lem3},~\ref{Lem4} and~\ref{Lem5}). By passing to the limit as the cutoff functions tend to the constant function equal to 1, we then obtain that both the torsion and Einstein curvature tensors of the contact form associated with the function $u$ (i.e. $u^{2/n}\Theta$, where $\Theta$ is the standard contact form on $\H^n$) vanish everywhere in $\H^n$, from which, as in~\cite{JL2}, we infer that $u$ must be of the form~\eqref{Th1Eq2}.

\section{Proofs of Theorems~\ref{Th1} and~\ref{Th2}}\label{Sec2}

This section is devoted to the proofs of Theorems~\ref{Th1} and~\ref{Th2}. Given a smooth function $u$ in $\H^n$, we denote $u_0=\partial_t u$ and for each $\alpha\in\left\{1,\dotsc,n\right\}$,
$$u_\alpha=u_{,\alpha}=Z_\alpha\(u\):=\partial_{z_\alpha}u+\sqrt{-1}\overline{z_\alpha}u_0$$
and
$$u_{\overline\alpha}=u_{,\overline\alpha}=Z_{\overline\alpha}\(u\):=\partial_{\overline{z_\alpha}}u-\sqrt{-1}z_\alpha u_0.$$
The vector fields $\left\{ Z_{\alpha},Z_{\bar \alpha}, \partial_{t} \right\}$ form a basis of left-invariant vector fields for the complex Lie algebra of $\mathbb{H}^{n}$, and they satisfy $\[Z_{\alpha},Z_{\bar\alpha}\] = -2\sqrt{-1}\partial_{t}$. We recall that in $\H^n$ equipped with its standard contact form, the Levi form is given by $2\delta_{\alpha\overline\beta}$, and the sub-Laplacian can be written as
$$\Delta_b u=-\Real u_{\alpha\overline\alpha}:=-\Real\sum_{\alpha=1}^nZ_{\overline\alpha}\(Z_\alpha\(u\)\).$$ 
We will write $\partial u$ to indicate the horizontal gradient of $u$ and 
$$B_{R}\(0\)=\left\{\(z,t\)\in\mathbb{H}^{n}:\,\(\left|z\right|^{4}+t^{2}\)^{1/4}<R\right\}$$
for the Kor\'anyi ball. Throughout this section, we also denote 
$$\left|V_\alpha\right|^2:=\sum_{\alpha=1}^nV_\alpha\overline{V_\alpha},\quad\left|S_{\alpha\beta}\right|^2:=\sum_{\alpha,\beta=1}^nS_{\alpha\beta}\overline{S_{\alpha\beta}}\quad\text{and}\quad\big|T_{\alpha\overline\beta}\big|^2:=\sum_{\alpha,\beta=1}^nT_{\alpha\overline\beta}\overline{T_{\alpha\overline\beta}}$$
for all tensors $V_\alpha$, $S_{\alpha\beta}$ and $T_{\alpha\overline\beta}$.

\smallskip
We let $u$ be a solution to \eqref{Th1Eq1}. We define
\begin{equation}\label{Th1Eq5}
f:=\frac{1}{n}\ln u-\ln2.
\end{equation}
It is easy to see that \eqref{Th1Eq1} can be rewritten as
\begin{equation}\label{Th1Eq6}
\Delta_bf=n\left|\partial f\right|^2+n\,e^{2f}\quad\text{in }\H^n.
\end{equation}
For each $\alpha,\beta\in\left\{1,\dotsc,n\right\}$, we define 
\begin{align*}
&D_{\alpha\beta}:=f_{\alpha\beta}-2f_\alpha f_\beta,\allowdisplaybreaks\\
&E_{\alpha\overline\beta}:=f_{\alpha\overline\beta}-\frac{1}{n}f_{\gamma\overline\gamma}\delta_{\alpha\overline\beta},\allowdisplaybreaks\\
&D_\alpha:=D_{\alpha\beta}f_{\overline\beta},\allowdisplaybreaks\\
&E_\alpha:=E_{\alpha\overline\beta}f_\beta\text{ and}\allowdisplaybreaks\\
&G_\alpha:=\sqrt{-1}f_{0\alpha}+gf_\alpha,
\end{align*}
where
$$g:=\left|\partial f\right|^2+e^{2f}-\sqrt{-1}f_0.$$
Up to some constant factors, $D_{\alpha\beta}$ and $E_{\alpha\overline\beta}$ correspond to the torsion and Einstein curvature tensors, respectively, of the contact form $u^{2/n}\Theta$, where $\Theta$ is the standard contact form on $\H^n$. We also observe that \eqref{Th1Eq6} can be rewritten as
\begin{equation}\label{Th1Eq7}
f_{\alpha\overline\alpha}+ng=0\quad\text{in }\H^n.
\end{equation}
We begin with proving a divergence formula. In the case where $m=0$, this formula is due to Jerison and Lee~\cite{JL2}. The idea of introducing the function $\left|g\right|^{-m}$ with $m>0$ in this formula is inspired from a similar and remarkable idea by Ou~\cite{Ou} in the context of the critical $p$--Laplace equation in $\R^n$. Our formula is as follows:

\begin{lemma}\label{Lem1}
Let $m\in\[0,1\]$, $n\in\N$ and $f$ be a real solution to \eqref{Th1Eq6}. Then
\begin{align}\label{Lem1Eq1}
&\Real\(\(e^{\(2n+m-2\)f}\left|g\right|^{-m}\(g\(D_\alpha+E_\alpha\)-\sqrt{-1}f_0\(D_\alpha-3E_\alpha+3G_\alpha\)\)\)_{\overline\alpha}\)\nonumber\\
&=e^{\(2n+m-2\)f}\left|g\right|^{-m}\Big(e^{2f}\(\big|D_{\alpha\beta}\big|^2+\big|E_{\alpha\overline\beta}\big|^2\)+\big|D_{\alpha\beta}f_{\overline\gamma}+E_{\alpha\overline\gamma}f_\beta\big|^2\nonumber\\
&\quad-\big|D_\alpha+E_\alpha\big|^2+c_1\left|G_\alpha+c_2D_\alpha+c_3E_\alpha\right|^2+c_4\left|D_\alpha+c_5E_\alpha\right|^2+c_6\left|E_\alpha\right|^2\Big),
\end{align}
where
\begin{align*}
c_1&:=3\left|g\right|^{-2}\(\left|g\right|^2-mf_0^2\),\allowdisplaybreaks\\
c_2&:=\frac{1}{3}\(1-\frac{2m\sqrt{-1}f_0\sqrt{\left|g\right|^2-f_0^2}}{\left|g\right|^2-mf_0^2}\),\allowdisplaybreaks\\
c_3&:=-\frac{1}{3}\(1+\frac{2m\sqrt{-1}f_0\sqrt{\left|g\right|^2-f_0^2}}{\left|g\right|^2-mf_0^2}\),\allowdisplaybreaks\\
c_4&:=\frac{1}{3}\(5-3m+\frac{4m\(1-m\)f_0^2}{\left|g\right|^2-mf_0^2}\),\allowdisplaybreaks\\
c_5&:=\frac{\(4-3m\)\left|g\right|^2-m\(2+m\)f_0^2+2m^2\left|g\right|^{-2}f_0^4}{\(5-3m\)\left|g\right|^2-m\(1+m\)f_0^2}\\
&\quad+\frac{2m\sqrt{-1}f_0\sqrt{\left|g\right|^2-f_0^2}\(\left|g\right|^2-mf_0^2\)}{\left|g\right|^2\(\(5-3m\)\left|g\right|^2-m\(1+m\)f_0^2\)}\text{ and}\allowdisplaybreaks\\
c_6&:=\frac{\(3-2m\)\left|g\right|^2+m\(5-6m\)f_0^2}{\(5-3m\)\left|g\right|^2-m\(1+m\)f_0^2}.
\end{align*}
\end{lemma}

Observe that since $0\le\left|f_0\right|<\left|g\right|$, we obtain that for each $m\in\[0,1\]$, $\left|g\right|^2-mf_0^2>0$ and $\(5-3m\)\left|g\right|^2-m\(1+m\)f_0^2>0$, hence $c_1$, $c_2$, $c_3$, $c_4$, $c_5$ and $c_6$ are well defined in this case.

\proof[Proof of Lemma~\ref{Lem1}]
Jerison and Lee's formula~\cite{JL2}*{Proposition~4.1} (see also~\cite{MO}*{Proposition~2.1}) gives
\begin{align}\label{Lem1Eq2}
&\Real\(\(e^{2\(n-1\)f}\(g\(D_\alpha+E_\alpha\)-\sqrt{-1}f_0\(D_\alpha-3E_\alpha+3G_\alpha\)\)\)_{\overline\alpha}\)\nonumber\\
&\quad=e^{2\(n-1\)f}\Big(e^{2f}\(\big|D_{\alpha\beta}\big|^2+\big|E_{\alpha\overline\beta}\big|^2\)+\left|G_\alpha\right|^2+\left|G_\alpha+D_\alpha\right|^2+\left|G_\alpha-E_\alpha\right|^2\nonumber\\
&\qquad+\big|D_{\alpha\beta}f_{\overline\gamma}+E_{\alpha\overline\gamma}f_\beta\big|^2\Big).
\end{align}
On the other hand, straightforward computations give
$$g_{\overline\alpha}=\overline{D_\alpha}+\overline{E_\alpha}+\overline{G_\alpha}\quad\text{and}\quad \overline{g}_{\overline\alpha}=\overline{D_\alpha}+\overline{E_\alpha}-\overline{G_\alpha}+2\overline{g}f_{\overline\alpha},$$
which in turn give
\begin{align}\label{Lem1Eq3}
\(e^{mf}\left|g\right|^{-m}\)_{\overline\alpha}&=me^{mf}\left|g\right|^{-m}\(f_{\overline\alpha}-\frac{1}{2}\left|g\right|^{-2}\(\overline{g}\,g_{\overline\alpha}+g\,\overline{g}_{\overline\alpha}\)\)\nonumber\\
&=-me^{mf}\left|g\right|^{-m-2}\Big(\(e^{2f}+\left|\partial f\right|^2\)\(\overline{D_\alpha}+\overline{E_\alpha}\)+\sqrt{-1}f_0\overline{G_\alpha}\Big).
\end{align}
It follows from \eqref{Lem1Eq2} and \eqref{Lem1Eq3} that
\begin{align}\label{Lem1Eq4}
&\Real\(\(e^{\(2n+m-2\)f}\left|g\right|^{-m}\(g\(D_\alpha+E_\alpha\)-\sqrt{-1}f_0\(D_\alpha-3E_\alpha+3G_\alpha\)\)\)_{\overline\alpha}\)\nonumber\\
&\quad=e^{\(2n+m-2\)f}\left|g\right|^{-m}\Big(e^{2f}\(\big|D_{\alpha\beta}\big|^2+\big|E_{\alpha\overline\beta}\big|^2\)+\big|D_{\alpha\beta}f_{\overline\gamma}+E_{\alpha\overline\gamma}f_\beta\big|^2\nonumber\\
&\qquad-\big|D_\alpha+E_\alpha\big|^2+\psi\Big).
\end{align}
where
\begin{align*}
\psi&:=\left|G_\alpha\right|^2+\left|G_\alpha+D_\alpha\right|^2+\left|G_\alpha-E_\alpha\right|^2+\left|D_\alpha+E_\alpha\right|^2\nonumber\\
&\quad+m\left|g\right|^{-2}\Real\(\sqrt{-1}f_0\(D_\alpha-3E_\alpha+3G_\alpha\)-g\(D_\alpha+E_\alpha\)\)\nonumber\\
&\qquad\times\Big(\(e^{2f}+\left|\partial f\right|^2\)\(\overline{D_\alpha}+\overline{E_\alpha}\)+\sqrt{-1}f_0\overline{G_\alpha}\Big)\nonumber\allowdisplaybreaks\\
&=\left|g\right|^{-2}\Big(3\(\left|g\right|^2-mf_0^2\)\left|G_\alpha\right|^2+\(\(2-m\)\left|g\right|^2+mf_0^2\)\(\left|D_\alpha\right|^2+\left|E_\alpha\right|^2\)\nonumber\allowdisplaybreaks\\
&\quad+2\(\left|g\right|^2-mf_0^2\)\Real\(G_\alpha\(\overline{D_\alpha}-\overline{E_\alpha}\)\)\nonumber\allowdisplaybreaks\\
&\quad-4mf_0\sqrt{\left|g\right|^2-f_0^2}\Imag\(G_\alpha\(\overline{D_\alpha}+\overline{E_\alpha}\)\)\nonumber\\
&\quad+2\(\(1-m\)\left|g\right|^2+mf_0^2\)\Real\(D_\alpha\overline{E_\alpha}\)-4mf_0\sqrt{\left|g\right|^2-f_0^2}\Imag\(D_\alpha\overline{E_\alpha}\)\Big).
\end{align*}
By completing the squares, we then obtain
\begin{equation}\label{Lem1Eq5}
\psi=c_1\left|G_\alpha+c_2D_\alpha+c_3E_\alpha\right|^2+c_4\left|D_\alpha+c_5E_\alpha\right|^2+c_6\left|E_\alpha\right|^2.
\end{equation}
Finally, by putting together \eqref{Lem1Eq4} and \eqref{Lem1Eq5}, we obtain \eqref{Lem1Eq1}.
\endproof

We now use Lemma~\ref{Lem1} to prove the following:

\begin{lemma}\label{Lem2}
Let $m\in\[0,1\)$, $n\in\N$, $f$ be a real solution to \eqref{Th1Eq6} and $\varphi$ be a smooth, nonnegative function with compact support in $\H^n$. Then
\begin{align}\label{Lem2Eq1}
\int_{\H^n}A_f\varphi&\le C\sqrt{\int_{\H^n}e^{\(2n+m-2\)f}\left|g\right|^{-m}\(e^{4f}+\left|\partial f\right|^4+f_0^2\)\varphi^{-1}\left|\partial\varphi\right|^2}\nonumber\\
&\quad\times\sqrt{\int_{\partial\varphi\ne0}A_f\varphi}
\end{align}
for some constant $C=C\(m\)>0$, where
\begin{equation}
A_f:=e^{\(2n+m-2\)f}\left|g\right|^{-m}\Big(e^{2f}\(\big|D_{\alpha\beta}\big|^2+\big|E_{\alpha\overline\beta}\big|^2\)+\left|D_\alpha\right|^2+\left|E_\alpha\right|^2+\left|G_\alpha\right|^2\Big).
\end{equation}
\end{lemma}

\proof[Proof of Lemma~\ref{Lem2}]
By multiplying \eqref{Lem1Eq1} by $\varphi$ and integrating by parts in $\H^n$, we obtain
\begin{align}\label{Lem2Eq2}
&\int_{\H^n}e^{\(2n+m-2\)f}\left|g\right|^{-m}\Big(e^{2f}\(\big|D_{\alpha\beta}\big|^2+\big|E_{\alpha\overline\beta}\big|^2\)+\big|D_{\alpha\beta}f_{\overline\gamma}+E_{\alpha\overline\gamma}f_\beta\big|^2\nonumber\\
&-\big|D_\alpha+E_\alpha\big|^2+c_1\left|G_\alpha+c_2D_\alpha+c_3E_\alpha\right|^2+c_4\left|D_\alpha+c_5E_\alpha\right|^2+c_6\left|E_\alpha\right|^2\Big)\varphi\nonumber\\
&=\Real\int_{\H^n}e^{\(2n+m-2\)f}\left|g\right|^{-m}\big(\sqrt{-1}f_0\(D_\alpha-3E_\alpha+3G_\alpha\)-g\(D_\alpha+E_\alpha\)\big)\varphi_{\overline\alpha},
\end{align}
where $c_1$, $c_2$, $c_3$, $c_4$, $c_5$ and $c_6$ are as in Lemma~\ref{Lem1}. We observe that 
\begin{align}\label{Lem2Eq3}
\big|D_{\alpha\beta}f_{\overline\gamma}+E_{\alpha\overline\gamma}f_\beta\big|^2&=\big|D_{\alpha\beta}f_{\overline\gamma}\big|^2+D_\alpha\overline{E_\alpha}+E_\alpha\overline{D_\alpha}+\big|E_{\alpha\overline\beta}f_\gamma\big|^2\nonumber\\
&=\big|D_{\alpha\beta}\big|^2\left|f_{\overline\gamma}\right|^2+D_\alpha\overline{E_\alpha}+E_\alpha\overline{D_\alpha}+\big|E_{\alpha\overline\beta}\big|^2\left|f_\gamma\right|^2\nonumber\\
&\ge\left|D_\alpha\right|^2+D_\alpha\overline{E_\alpha}+E_\alpha\overline{D_\alpha}+\left|E_\alpha\right|^2\nonumber\\
&=\left|D_\alpha+E_\alpha\right|^2.
\end{align} 
Moreover, since $0\le m<1$ and $0\le\left|f_0\right|<\left|g\right|$, straightforward computations give
\begin{align}
&c_1\ge3\(1-m\)>0,\allowdisplaybreaks\label{Lem2Eq4}\\
&c_4\ge\frac{5-3m}{3}>0,\allowdisplaybreaks\label{Lem2Eq5}\\
&c_6\ge\frac{3-2m}{5-3m}>0\text{ and}\allowdisplaybreaks\label{Lem2Eq6}\\
&\left|c_2\right|+\left|c_3\right|+\left|c_5\right|\le C\label{Lem2Eq7}
\end{align}
for some constant $C=C\(m\)>0$. It follows from \eqref{Lem2Eq4}, \eqref{Lem2Eq5}, \eqref{Lem2Eq6} and \eqref{Lem2Eq7} that 
\begin{multline}\label{Lem2Eq8}
c_1\left|G_\alpha+c_2D_\alpha+c_3E_\alpha\right|^2+c_4\left|D_\alpha+c_5E_\alpha\right|^2+c_6\left|E_\alpha\right|^2\\
\ge C\(\left|D_\alpha\right|^2+\left|E_\alpha\right|^2+\left|G_\alpha\right|^2\).
\end{multline}
for some constant $C=C\(m\)>0$. On the other hand, straightforward estimates together with Cauchy--Schwartz' and Young's inequalities give 
\begin{align}\label{Lem2Eq9}
&\Real\int_{\H^n}e^{\(2n+m-2\)f}\left|g\right|^{-m}\(\sqrt{-1}f_0\(D_\alpha-3E_\alpha+3G_\alpha\)-g\(D_\alpha+E_\alpha\)\)\varphi_{\overline\alpha}\nonumber\\
&\quad\le\int_{\H^n}e^{\(2n+m-2\)f}\left|g\right|^{-m}\Big(\(e^{2f}+\left|\partial f\right|^2\)\(\left|D_\alpha\right|+\left|E_\alpha\right|\)\nonumber\\
&\qquad+\left|f_0\right|\(2\left|D_\alpha\right|+2\left|E_\alpha\right|+3\left|G_\alpha\right|\)\Big)\left|\partial\varphi\right|\nonumber\allowdisplaybreaks\\
&\quad\le\int_{\H^n}e^{\(2n+m-2\)f}\left|g\right|^{-m}\Big(e^{2f}\left|\partial f\right|\(\big|D_{\alpha\beta}\big|+\big|E_{\alpha\overline\beta}\big|\)+\left|\partial f\right|^2\(\left|D_\alpha\right|+\left|E_\alpha\right|\)\nonumber\\
&\qquad+\left|f_0\right|\(2\left|D_\alpha\right|+2\left|E_\alpha\right|+3\left|G_\alpha\right|\)\Big)\left|\partial\varphi\right|\nonumber\allowdisplaybreaks\\
&\quad\le\frac{1}{2}\int_{\H^n}e^{\(2n+m-2\)f}\left|g\right|^{-m}\Big(\(e^{2f}+\left|\partial f\right|^2\)e^f\(\big|D_{\alpha\beta}\big|+\big|E_{\alpha\overline\beta}\big|\)\nonumber\\
&\qquad+2\left|\partial f\right|^2\(\left|D_\alpha\right|+\left|E_\alpha\right|\)+2\left|f_0\right|\(2\left|D_\alpha\right|+2\left|E_\alpha\right|+3\left|G_\alpha\right|\)\Big)\left|\partial\varphi\right|\nonumber\allowdisplaybreaks\\
&\quad\le11\sqrt{\int_{\H^n}e^{\(2n+m-2\)f}\left|g\right|^{-m}\(e^{4f}+\left|\partial f\right|^4+f_0^2\)\varphi^{-1}\left|\partial\varphi\right|^2}\nonumber\\
&\qquad\times\sqrt{\int_{\partial\varphi\ne0}A_f\varphi}.
\end{align}
Finally, by combining \eqref{Lem2Eq2}, \eqref{Lem2Eq3}, \eqref{Lem2Eq8} and \eqref{Lem2Eq9}, we obtain \eqref{Lem2Eq1}.
\endproof

The next three lemmas will be used to estimates the terms in the right-hand side of \eqref{Lem2Eq1}. First, we prove the following:

\begin{lemma}\label{Lem3}
Let $m\in\[0,1\]$, $\sigma\in\R$, $\varepsilon>0$, $n\in\N$, $f$ be a real solution to \eqref{Th1Eq5} and $\varphi$ be a smooth, nonnegative function with compact support in $\H^n$. Then 
\begin{align}\label{Lem3Eq1}
&\int_{\H^n}e^{\(2n+m-2\)f}\left|g\right|^{-m}f_0^2\varphi\nonumber\\
&\quad\le C\int_{\H^n}e^{\(2n+m-2\)f}\left|g\right|^{-m}\Big(e^{4f}\varphi+\left|\partial f\right|^4\varphi+\varphi^{-3}\left|\partial\varphi\right|^4+\varepsilon^{-2}\varphi^{1-2\sigma}\nonumber\\
&\qquad+\varepsilon\(\left|D_\alpha\right|^2+\left|E_\alpha\right|^2+\left|G_\alpha\right|^2\)\varphi^{1+\sigma}\Big)
\end{align}
for some constant $C=C\(n\)>0$.
\end{lemma}

\proof[Proof of Lemma~\ref{Lem3}]
By observing that $f_{0\overline\alpha}=\sqrt{-1}\,\overline{G_\alpha}-\sqrt{-1}\,\overline{g}f_{\overline\alpha}$ and using \eqref{Th1Eq7} and \eqref{Lem1Eq3}, we obtain
\begin{align}\label{Lem3Eq2}
&\Imag\(\(e^{\(2n+m-2\)f}\left|g\right|^{-m}f_0f_\alpha\varphi\)_{\overline\alpha}\)\nonumber\\
&\quad=e^{\(2n+m-2\)f}\left|g\right|^{-m}\(nf_0^2-e^{2f}\left|\partial f\right|^2-\left|\partial f\right|^4\)\varphi\nonumber\allowdisplaybreaks\\
&\qquad+\Imag\Big(e^{\(2n+m-2\)f}\left|g\right|^{-m}f_\alpha\big(\left|g\right|^{-2}\big(\sqrt{-1}\(\left|g\right|^2-mf_0^2\)\overline{G_\alpha}\nonumber\\
&\qquad-mf_0\(e^{2f}+\left|\partial f\right|^2\)\(\overline{D_\alpha}+\overline{E_\alpha}\)\big)\varphi+f_0\varphi_{\overline\alpha}\big)\Big).
\end{align}
By integrating \eqref{Lem3Eq2} in $\H^n$, we obtain
\begin{align}\label{Lem3Eq3}
&\int_{\H^n}e^{\(2n+m-2\)f}\left|g\right|^{-m}f_0^2\varphi\nonumber\\
&\quad=\frac{1}{n}\bigg(\int_{\H^n}e^{\(2n+m-2\)f}\left|g\right|^{-m}\(e^{2f}\left|\partial f\right|^2+\left|\partial f\right|^4\)\varphi\nonumber\allowdisplaybreaks\\
&\qquad-\Imag\int_{\H^n}e^{\(2n+m-2\)f}\left|g\right|^{-m}f_\alpha\Big(\left|g\right|^{-2}\Big(\sqrt{-1}\(\left|g\right|^2-mf_0^2\)\overline{G_\alpha}\nonumber\\
&\qquad-mf_0\(e^{2f}+\left|\partial f\right|^2\)\(\overline{D_\alpha}+\overline{E_\alpha}\)\Big)\varphi+f_0\varphi_{\overline\alpha}\Big)\bigg).
\end{align}
Straightforward estimates together with Young's inequality give 
\begin{align}\label{Lem3Eq4}
\int_{\H^n}e^{\(2n+m\)f}\left|g\right|^{-m}\left|\partial f\right|^2\varphi\le\frac{1}{2}\int_{\H^n}e^{\(2n+m-2\)f}\left|g\right|^{-m}\(e^{4f}+\left|\partial f\right|^4\)\varphi
\end{align}
and
\begin{align}\label{Lem3Eq5}
&-\Imag\int_{\H^n}e^{\(2n+m-2\)f}\left|g\right|^{-m}f_0f_\alpha\varphi_{\overline\alpha}\nonumber\\
&\quad\le\int_{\H^n}e^{\(2n+m-2\)f}\left|g\right|^{-m}\left|f_0\right|\left|\partial f\right|\left|\partial\varphi\right|\nonumber\allowdisplaybreaks\\
&\quad\le\frac{1}{4}\int_{\H^n}e^{\(2n+m-2\)f}\left|g\right|^{-m}\(2f_0^2\varphi+\left|\partial f\right|^4\varphi+\varphi^{-3}\left|\partial\varphi\right|^4\).
\end{align}
In a similar way, we obtain
\begin{align}\label{Lem3Eq6}
&-\Imag\int_{\H^n}e^{\(2n+m-2\)f}\left|g\right|^{-m-2}f_\alpha\Big(\sqrt{-1}\(\left|g\right|^2-mf_0^2\)\overline{G_\alpha}\nonumber\\
&\qquad-mf_0\(e^{2f}+\left|\partial f\right|^2\)\(\overline{D_\alpha}+\overline{E_\alpha}\)\Big)\varphi\nonumber\\
&\quad\le \int_{\H^n}e^{\(2n+m-2\)f}\left|g\right|^{-m}\left|\partial f\right|\(\left|D_\alpha\right|+\left|E_\alpha\right|+\left|G_\alpha\right|\)\varphi\nonumber\allowdisplaybreaks\\
&\quad\le\frac{1}{4}\int_{\H^n}e^{\(2n+m-2\)f}\left|g\right|^{-m}\Big(3\left|\partial f\right|^4\varphi+3\varepsilon^{-2}\varphi^{1-2\sigma}\nonumber\\
&\qquad+2\varepsilon\(\left|D_\alpha\right|^2+\left|E_\alpha\right|^2+\left|G_\alpha\right|^2\)\varphi^{1+\sigma}\Big).
\end{align}
Finally, by combining \eqref{Lem3Eq3}, \eqref{Lem3Eq4}, \eqref{Lem3Eq5} and \eqref{Lem3Eq6}, we obtain \eqref{Lem3Eq1}.
\endproof

Now, we prove the following:

\begin{lemma}\label{Lem4}
Let $m\in\[0,1\]$, $\sigma\in\R$, $\varepsilon>0$, $n\ge2$, $f$ be a real solution to \eqref{Th1Eq5} and $\varphi$ be a smooth, nonnegative function with compact support in $\H^n$. Then 
\begin{align}\label{Lem4Eq1}
&\int_{\H^n}e^{\(2n+m-2\)f}\left|g\right|^{-m}\left|\partial f\right|^4\varphi\nonumber\\
&\quad\le C\int_{\H^n}e^{\(2n+m-2\)f}\left|g\right|^{-m}\Big(e^{4f}\varphi+\varphi^{-3}\left|\partial\varphi\right|^4+\varepsilon^{-2}\varphi^{1-2\sigma}\nonumber\\
&\qquad+\varepsilon\(\left|D_\alpha\right|^2+\left|E_\alpha\right|^2+\left|G_\alpha\right|^2\)\varphi^{1+\sigma}\Big)
\end{align}
for some constant $C=C\(n\)>0$.
\end{lemma}

\proof[Proof of Lemma~\ref{Lem4}]
By using \eqref{Th1Eq7} and \eqref{Lem1Eq3} together with $f_{\overline\alpha\overline\beta}f_\alpha f_\beta=f_\alpha\overline{D_\alpha}+2\left|\partial f\right|^4$ and $f_{\beta\overline\alpha}f_\alpha f_{\overline\beta}=f_\alpha\overline{E_\alpha}-g\left|\partial f\right|^2$, we obtain
\begin{align}\label{Lem4Eq2}
&\Real\(\(e^{\(2n+m-2\)f}\left|g\right|^{-m}\left|\partial f\right|^2f_\alpha\varphi\)_{\overline\alpha}\)\nonumber\\
&\quad=e^{\(2n+m-2\)f}\left|g\right|^{-m}\(\(n-1\)\left|\partial f\right|^4-\(n+1\)e^{2f}\left|\partial f\right|^2\)\varphi\nonumber\allowdisplaybreaks\\
&\qquad+\Real\Big(e^{\(2n+m-2\)f}\left|g\right|^{-m} f_\alpha\Big(\left|g\right|^{-2}\Big(\(\left|g\right|^2-me^{2f}\left|\partial f\right|^2-m\left|\partial f\right|^4\)\nonumber\\
&\qquad\times\(\overline{D_\alpha}+\overline{E_\alpha}\)-m\sqrt{-1}\left|\partial f\right|^2f_0\overline{G_\alpha}\Big)\varphi+\left|\partial f\right|^2\varphi_{\overline\alpha}\Big)\Big).
\end{align}
By integrating \eqref{Lem4Eq2} in $\H^n$, we obtain
\begin{align}\label{Lem4Eq3}
&\int_{\H^n}e^{\(2n+m-2\)f}\left|g\right|^{-m}\left|\partial f\right|^4\varphi\nonumber\\
&\quad=\frac{1}{n-1}\bigg(\(n+1\)\int_{\H^n}e^{\(2n+m\)f}\left|g\right|^{-m}\left|\partial f\right|^2\varphi\nonumber\allowdisplaybreaks\\
&\qquad-\Real\int_{\H^n}e^{\(2n+m-2\)f}\left|g\right|^{-m}f_\alpha\Big(\left|g\right|^{-2}\Big(\Big(\left|g\right|^2-me^{2f}\left|\partial f\right|^2-m\left|\partial f\right|^4\Big)\nonumber\\
&\qquad\times\(\overline{D_\alpha}+\overline{E_\alpha}\)-m\sqrt{-1}\left|\partial f\right|^2f_0\overline{G_\alpha}\Big)\varphi+\left|\partial f\right|^2\varphi_{\overline\alpha}\Big)\bigg).
\end{align}
For each $\rho>0$, straightforward estimates together with Young's inequality give 
\begin{align}\label{Lem4Eq4}
\int_{\H^n}e^{\(2n+m\)f}\left|g\right|^{-m}\left|\partial f\right|^2\varphi\le\frac{1}{2}\int_{\H^n}e^{\(2n+m-2\)f}\left|g\right|^{-m}\(\rho^{-1}e^{4f}+\rho\left|\partial f\right|^4\)\varphi
\end{align}
and
\begin{align}\label{Lem4Eq5}
&-\Real\int_{\H^n}e^{\(2n+m-2\)f}\left|g\right|^{-m}\left|\partial f\right|^2f_\alpha\varphi_{\overline\alpha}\nonumber\\
&\quad\le\int_{\H^n}e^{\(2n+m-2\)f}\left|g\right|^{-m}\left|\partial f\right|^3\left|\partial\varphi\right|\nonumber\\
&\quad\le\frac{1}{4}\int_{\H^n}e^{\(2n+m-2\)f}\left|g\right|^{-m}\(3\rho\left|\partial f\right|^4\varphi+\rho^{-3}\varphi^{-3}\left|\partial\varphi\right|^4\).
\end{align}
In a similar way, we obtain
\begin{align}\label{Lem4Eq6}
&-\Real\int_{\H^n}e^{\(2n+m-2\)f}\left|g\right|^{-m-2}f_\alpha\Big(\Big(\left|g\right|^2-me^{2f}\left|\partial f\right|^2-m\left|\partial f\right|^4\Big)\(\overline{D_\alpha}+\overline{E_\alpha}\)\nonumber\\
&\qquad-m\sqrt{-1}\left|\partial f\right|^2f_0\overline{G_\alpha}\Big)\varphi\Big)\nonumber\\
&\quad\le\int_{\H^n}e^{\(2n+m-2\)f}\left|g\right|^{-m}\left|\partial f\right|\(\left|D_\alpha\right|+\left|E_\alpha\right|+\left|G_\alpha\right|\)\varphi\nonumber\allowdisplaybreaks\\
&\quad\le\frac{1}{4}\int_{\H^n}e^{\(2n+m-2\)f}\left|g\right|^{-m}\Big(3\rho\left|\partial f\right|^4\varphi+3\rho^{-1}\varepsilon^{-2}\varphi^{1-2\sigma}\nonumber\\
&\qquad+2\varepsilon\(\left|D_\alpha\right|^2+\left|E_\alpha\right|^2+\left|G_\alpha\right|^2\)\varphi^{1+\sigma}\Big).
\end{align}
Finally, by putting together \eqref{Lem4Eq3}, \eqref{Lem4Eq4}, \eqref{Lem4Eq5} and \eqref{Lem4Eq6}, and letting $\rho$ be small enough, we obtain \eqref{Lem4Eq1}.
\endproof

Next, we prove the following:

\begin{lemma}\label{Lem5}
Let $R>1$, $n\in\N$, $r\in\[0,2\]$, $q\in\[0,n+2\]$ if $r=0$, $q\in\[0,n+2-r\)$ if $r\in\(0,2\]$ and $f$ be a real solution to \eqref{Th1Eq6}. Then
\begin{equation}\label{Lem5Eq1}
\int_{B_R\(0\)}e^{qf}\left|\partial f\right|^r\le CR^{2n+2-q-r},
\end{equation}
for some constant $C=C\(n,q,r\)>0$.
\end{lemma}

\proof[Proof of Lemma~\ref{Lem5}]
We let $\eta$ be a smooth cutoff function in $\H^n$ such that $\eta\equiv1$ in $B_1\(0\)$, $\eta\equiv0$ in $\H^n\backslash B_2\(0\)$ and $0\le\eta\le1$ in $B_2\(0\)\backslash B_1\(0\)$. For each $R>1$, we define
$$\eta_R\(z,t\):=\eta\(R^{-1}z,R^{-2}t\)\quad\forall\(z,t\)\in\H^n.$$
For each $\theta>1$, by using \eqref{Th1Eq7}, we obtain
\begin{equation}\label{Lem5Eq2}
\Real\(\(e^{qf}f_\alpha\eta_R^\theta\)_{\overline\alpha}\)=e^{qf}\(\(q-n\)\left|\partial f\right|^2-ne^{2f}\)\eta_R^\theta+\theta e^{qf}\eta_R^{\theta-1}f_\alpha\(\eta_R\)_{\overline\alpha}.
\end{equation}
By integrating \eqref{Lem5Eq2} in $\H^n$, we obtain
\begin{equation}\label{Lem5Eq3}
\int_{\H^n}e^{qf}\(\(n-q\)\left|\partial f\right|^2+ne^{2f}\)\eta_R^\theta=\theta\int_{\H^n}e^{qf}\eta_R^{\theta-1}f_\alpha\(\eta_R\)_{\overline\alpha}.
\end{equation}
For each $\rho>0$, straightforward estimates together with Young's inequality give 
\begin{align}
\int_{\H^n}e^{qf}\eta_R^{\theta-1}f_\alpha\(\eta_R\)_{\overline\alpha}&\le\int_{\H^n}e^{qf}\left|\partial f\right|\eta_R^{\theta-1}\left|\partial\eta_R\right|\label{Lem5Eq4}\\
&\le \frac{1}{4\rho}\int_{\H^n}e^{qf}\eta_R^{\theta-2}\left|\partial\eta_R\right|^2+\rho\int_{\H^n}e^{qf}\left|\partial f\right|^2\eta_R^\theta\label{Lem5Eq5}
\end{align}
provided we choose $\theta>2$. If $q>0$, then another application of Young's inequality gives
\begin{equation}\label{Lem5Eq6}
\int_{\H^n}e^{qf}\eta_R^{\theta-2}\left|\partial\eta_R\right|^2\le C_\rho\int_{\H^n}\eta_R^{\theta-q-2}\left|\partial\eta_R\right|^{q+2}+\rho\int_{\H^n}e^{\(q+2\)f}\eta_R^\theta
\end{equation}
for some constant $C_\rho=C_\rho\(q\)>0$, provided we choose $\theta>q+2$. Moreover, straightforward estimates give
\begin{equation}\label{Lem5Eq7}
\int_{\H^n}\eta_R^{\theta-q-2}\left|\partial\eta_R\right|^{q+2}\le C R^{2n-q}
\end{equation}
for some constant $C=C\(n,q\)>0$. By choosing $\rho$ small enough and putting together \eqref{Lem5Eq3}, \eqref{Lem5Eq5}, \eqref{Lem5Eq6} and \eqref{Lem5Eq7}, we obtain \eqref{Lem5Eq1} for $q\in\[0,n\)$ when $r=2$ and $q\in\[2,n+2\)$ when $r=0$. In the case where $0<r<2$ and $2-r\le q<n+2-r$, H\"older's inequality gives 
\begin{align*}
\int_{B_R\(0\)}e^{qf}\left|\partial f\right|^r&\le\(\int_{B_R\(0\)}e^{\(q+r\)f}\)^{\frac{2-r}{2}}\(\int_{B_R\(0\)}e^{\(q+r-2\)f}\left|\partial f\right|^2\)^{\frac{r}{2}}\\
&\le C\(R^{2n+2-q-r}\)^{\frac{2-r}{2}}\(R^{2n+2-q-r}\)^{\frac{r}{2}}\\
&=CR^{2n+2-q-r}
\end{align*}
for some constant $C=C\(n,q,r\)>0$, thus \eqref{Lem5Eq1} still holds in this case. In the case where $0\le r<2$ and $0\le q<2-r$, H\"older's inequality gives 
\begin{align*}
\int_{B_R\(0\)}e^{qf}\left|\partial f\right|^r&\le\(\int_{B_R\(0\)}1\)^{\frac{2-q-r}{2}}\(\int_{B_R\(0\)}e^{\frac{2q}{q+r}f}\left|\partial f\right|^{\frac{2r}{q+r}}\)^{\frac{q+r}{2}}\\
&\le C\(R^{2n+2}\)^{\frac{2-q-r}{2}}\(R^{2n}\)^{\frac{q+r}{2}}\\
&=CR^{2n+2-q-r}
\end{align*}
for some constant $C=C\(n,q,r\)>0$, thus \eqref{Lem5Eq1} still holds in this case. Finally, in the case where $r=0$ and $q=n+2$, by using \eqref{Lem5Eq3} and \eqref{Lem5Eq4}, we obtain
$$\int_{B_R\(0\)}e^{\(n+2\)f}\le CR^{-1}\int_{B_{2R}\(0\)}e^{nf}\left|\partial f\right|\le C'R^n$$
for some constants $C=C\(n\)>0$ and $C'=C'\(n\)>0$, thus \eqref{Lem5Eq1} still holds in this case.
\endproof

Now, we can prove Theorems~\ref{Th1} and~\ref{Th2}, starting with Theorem~\ref{Th2}:

\proof[Proof of Theorem~\ref{Th2}]
We let $u$ be a solution to \eqref{Th1Eq1} satisfying \eqref{Th2Eq1}. We let $f$ be as in \eqref{Th1Eq5}. We let $\eta$ and $\eta_R$ be as in the proof of Lemma~\ref{Lem5}. For each $m\in\[0,1\)$ and $\theta>2$, by applying Lemma~\ref{Lem2} with $\varphi=\eta_R^\theta$, we obtain
\begin{align}\label{Th1Eq8}
\int_{\H^n}A_f\eta_R^\theta&\le CR^{-1}\sqrt{\int_{\H^n}e^{\(2n+m-2\)f}\left|g\right|^{-m}\(e^{4f}+\left|\partial f\right|^4+f_0^2\)\eta_R^{\theta-2}}\nonumber\\
&\quad\times\sqrt{\int_{B_{2R}\(0\)\backslash B_R\(0\)}A_f\eta_R^\theta}
\end{align}
for some constant $C=C\(m,\theta,\left\|\partial\eta\right\|_\infty\)>0$. By using \eqref{Th1Eq8} and Lemma~\ref{Lem3} with $\varphi=\eta_R^{\theta-2}$, $\sigma=\frac{2}{\theta-2}$ and $\varepsilon=\rho R^2$ for some small $\rho>0$, we obtain
\begin{align}\label{Th1Eq9}
\(\int_{\H^n}A_f\eta_R^\theta\)^2&\le CR^{-2}\bigg(\int_{\H^n}e^{\(2n+m-2\)f}\left|g\right|^{-m}\Big(e^{4f}\eta_R^{\theta-2}+\left|\partial f\right|^4\eta_R^{\theta-2}\nonumber\\
&\quad+R^{-4}\eta_R^{\theta-4}\Big)\bigg)\(\int_{B_{2R}\(0\)\backslash B_R\(0\)}A_f\eta_R^\theta\)
\end{align}
for some constant $C=C\(n,m,\theta,\left\|\partial\eta\right\|_\infty\)>0$, provided we choose $\theta>4$. By using \eqref{Th1Eq9} and Lemma~\ref{Lem4} with $\varphi=\eta_R^{\theta-2}$, $\sigma=\frac{2}{\theta-2}$ and $\varepsilon=\rho R^2$ for some small $\rho>0$, we obtain
\begin{align}\label{Th1Eq10}
\(\int_{\H^n}A_f\eta_R^\theta\)^2&\le CR^{-2}\(\int_{\H^n}e^{\(2n+m+2\)f}\left|g\right|^{-m}\(\eta_R^{\theta-2}+R^{-4}e^{-4f}\eta_R^{\theta-4}\)\)\nonumber\\
&\quad\times\(\int_{B_{2R}\(0\)\backslash B_R\(0\)}A_f\eta_R^\theta\)\nonumber\allowdisplaybreaks\\
&\le CR^{-2}\(\int_{B_{2R}\(0\)}e^{\(2n-m+2\)f}\(1+R^{-4}e^{-4f}\)\)\nonumber\\
&\quad\times\(\int_{B_{2R}\(0\)\backslash B_R\(0\)}A_f\eta_R^\theta\)
\end{align}
for some constant $C=C\(n,m,\theta,\left\|\partial\eta\right\|_\infty\)>0$. By letting $m=n\(2-q\)+2$ and using \eqref{Th2Eq1}, we obtain
\begin{equation}\label{Th1Eq11}
\int_{B_{2R}\(0\)}e^{\(2n-m+2\)f}\le CR^2
\end{equation}
for some constant $C=C\(n,q,f\)>0$. Moreover, in the case where $2n-m-2>n+2$, i.e. $n-m>4$, by using H\"older's inequality together with \eqref{Th1Eq11} and Lemma~\ref{Lem5}, we obtain
\begin{align}\label{Th1Eq12}
\int_{B_{2R}\(0\)}e^{\(2n-m-2\)f}&\le\(\int_{B_{2R}\(0\)}e^{\(2n-m+2\)f}\)^{\frac{n-m-4}{n-m}}\(\int_{B_{2R}\(0\)}e^{\(n+2\)f}\)^{\frac{4}{n-m}}\nonumber\\
&\le C\(R^2\)^{\frac{n-m-4}{n-m}}\(R^n\)^{\frac{4}{n-m}}\nonumber\\
&=CR^{6-\frac{4\(2-m\)}{n-m}}\nonumber\\
&=\smallo\(R^6\)\quad\text{as }R\to\infty.
\end{align}
When $2n-m-2\le n+2$, a direct application of Lemma~\ref{Lem5} gives
\begin{align}\label{Th1Eq13}
\int_{B_{2R}\(0\)}e^{\(2n-m-2\)f}&\le CR^{4+m}\nonumber\\
&=\smallo\(R^6\)\quad\text{as }R\to\infty.
\end{align}
It follows from \eqref{Th1Eq10}, \eqref{Th1Eq11}, \eqref{Th1Eq12} and \eqref{Th1Eq13} that
$$\(\int_{\H^n}A_f\eta_R^\theta\)^2\le C\int_{B_{2R}\(0\)\backslash B_R\(0\)}A_f\eta_R^\theta$$
for some constant $C=C\(n,q,\theta,\rho,f\)>0$, which in turn gives
\begin{equation}\label{Th1Eq14}
\int_{\H^n}A_f\eta_R^\theta\to\int_{\H^n}A_f=0\quad\text{as }R\to\infty.
\end{equation}
It follows from \eqref{Th1Eq14} that $D_{\alpha\beta}\equiv0$ and $E_{\alpha\overline\beta}\equiv0$ in $\H^n$ for all $\alpha,\beta\in\left\{1,2\right\}$. We then conclude as in~\cite{JL2} that $u$ is of the form
\endproof

\proof[Proof of Theorem~\ref{Th1}]
By using \eqref{Th1Eq3} and Lemma~\ref{Lem5}, we obtain 
\begin{align*}
\int_{B_{2R}\(0\)}u^{\frac{2n+2}{n}}&\le CR^{2-n}\int_{B_{2R}\(0\)}u^{\frac{n+2}{n}}\\
&\le C'R^2
\end{align*}
for some constants $C=C\(n,u\)>0$ and $C'=C'\(n,u\)>0$. By applying Theorem~\ref{Th2}, we then obtain that $u$ is of the form \eqref{Th1Eq2}.
\endproof

\end{document}